\documentclass[reqno]{amsart}
 \usepackage{amsbsy,amssymb,amscd,amsfonts,latexsym,amstext,delarray,
 amsmath,graphicx, color}
 \usepackage[dvips]{epsfig}
\usepackage{hyperref}
\input xypic

\definecolor{Green}{rgb}{0.0,0.40,0.0}

\newtheorem{thm}{Theorem}[section]
\newtheorem{prop}[thm]{Proposition}

\newtheorem{lem}[thm]{Lemma}

\newtheorem{defn}[thm]{Definition}

\def\qqq{\,,\quad~\forall}

\def\GL{{\rm GL}}

\def\Hom{{\rm Hom}}

\def\Spec{{\rm Spec}}
\def\Sp{{\rm Spec}}

\def\A{{\mathbb A}}
\def\C{{\mathbb C}}
\def\F{{\mathbb F}}
\def\K{{\mathbb K}}

\def\N{{\mathbb N}}
\def\P{{\mathbb P}}
\def\Q{{\mathbb Q}}
\def\R{{\mathbb R}}
\def\Z{{\mathbb Z}}

\def\cU{{\cal U}}

\def\cA{{\mathcal A}}

\def\cG{{\mathcal G}}

\def\cF{{\mathcal F}}

\def\cM{{\mathcal M}}
\def\cN{{\mathcal N}}
\def\cO{{\mathcal O}}

\def\cR{{\mathcal R}}

\def\cT{{\mathcal T}}

\newcommand{\ie}{{\it i.e.\/}\ }
\newcommand{\eg}{{\it e.g.\/}\ }
\newcommand{\cf}{{\it cf.\/}\ }
\newcommand{\opcit}{{\it op.cit.\/}\ }
\newcommand{\resp}{{\it resp.\/}\ }

\def\cal{\mathcal}

\def\dim{{\mbox{dim}}}

\def\Hom {{\mbox{Hom}}}

\def\bG{\mathbb G}
\def\Gm{{\mathbb G}_m}

\def\abgr{D}
\def\pp{{\rm prod}}

\newcommand{\nil}[1]{}

\title
{On the notion of geometry over $\F_1$}

\author[Connes]{Alain Connes}
\author[Consani]{Caterina Consani}
\address{A.~Connes: Coll\`ege de France \\
3, rue d'Ulm \\ Paris, F-75005 France
\\ I.H.E.S. and Vanderbilt
University} \email{alain\@@connes.org}
\address{C.~Consani: Mathematics Department \\ Johns Hopkins
University \\ Baltimore, MD 21218 USA} \email{kc\@@math.jhu.edu}

\date{}
\begin{document}
\maketitle \vspace{2cm}

\begin{abstract} We refine the notion
of variety over the ``field with one element" developed by C. Soul\'e by
introducing a grading in the associated functor to the category of sets, and
show that this notion becomes compatible with the geometric viewpoint developed
by J. Tits. We then solve an open question of C. Soul\'e by proving, using
results of J. Tits and C. Chevalley, that Chevalley group schemes are examples of
varieties over a quadratic extension of the above ``field".
\end{abstract}

 \tableofcontents

  \section{Introduction}

In his theory of buildings  J.~Tits obtained a broad generalization of the
celebrated von Staudt reconstruction theorem in projective geometry, involving
as groups of symmetries not only $\GL_n$ but the full collection of Chevalley
algebraic groups. Among the axioms (\cite{Tits2}) which characterize these constructions, a
relevant one is played by the condition of ``thickness" which  states, in its
simplest form, that a projective line contains at least three points. By
replacing this requirement with its strong negation, \ie by imposing that a
line contains exactly two points, one still obtains a coherent ``geometry"
which is a degenerate form of classical projective geometry. In the case of
buildings, this degenerate case is described by the theory of ``thin" complexes
and in particular by the structure of the apartments, which are the basic
constituents of the theory of buildings. The degeneracy of the von Staudt field
inspired to Tits the conviction that these degenerate forms of geometries are a
manifestation of the existence of a hypothetical algebraic object that he named
``the field of characteristic one'' (\cite{Tits}). The richness and beauty of this geometric
picture gives convincing evidence for the pertinence of a separate study of the
degenerate case.

\smallskip

For completely independent reasons the need for  a ``field of characteristic
one" (also called the field with one element) has also emerged in Arakelov
geometry and \eg in \cite{Manin}, in the context of a geometric interpretation
of the zeros of zeta and $L$-functions. These speculative constructions aim
for the description of  a geometric framework analogous to the one used by Weil
in the proof of the Riemann Hypothesis for function fields. More precisely, one
seeks for a replacement of the surface $C\times_{\F_q} C$, where $C$ is a
(projective, smooth algebraic) curve over a finite field $\F_q$ and whose field
of functions is the given global field. The main idea  is   to postulate the
existence of the ``absolute point" $\Spec\,\F_1$ over which  any algebraic
scheme would sit. In the particular case of $\Spec\,\Z$, one would then be able
to use the spectrum of the tensor product $\Z \times_{\F_1}\Z$ as a substitute
for the surface $C\times_{\F_q} C$. This viewpoint has given rise, in the
recent past, to a number of interesting constructions (\cf \cite{Kapranov},
\cite{Manin}, \cite{Soule99}, \cite{Soule}, \cite{Durov}, \cite{Haran},
\cite{TV}).

\smallskip

Our interest in the quest for $\F_1$  arose from the following equation
\begin{equation}\label{basic}
\F_{1^n}\otimes_{\F_1} \Z:=\Z[T]/(T^n-1)\,, \qquad n\in\N
\end{equation}
which was introduced in \cite{Soule} and supplies a definition of the finite
extension $\F_{1^n}$ of $\F_1$, after base change to $\Z$. The main point
promoted in \cite{ccm} is that the above equation \eqref{basic} yields (without
knowing the algebraic nature of $\F_1$, and after base change to $\Z$), an
algebraic object which reflects the structure of the inductive limit
$\F_{1^\infty} = \varinjlim \F_{1^n} $, by supplying also an analogue of the
geometric Frobenius correspondence. This object is the integral model of a
rational Hecke algebra  which defines  the quantum statistical mechanical
system of  \cite{BC}. It is known that our construction determines, after
passing to the dual system, a spectral realization of the zeros of the Riemann
zeta function, as well as a trace formula interpretation of the Riemann-Weil
explicit formulas (see \cite{BC}, \cite{Co-zeta}, \cite{CCM}, \cite{CCM2},
\cite{Meyer}).

\smallskip
In \cite{ccm}, we made use of the general definition of an algebraic variety
over $\F_1$  as introduced by C. Soul\'e in \cite{Soule}. Our goal in this
paper is to give an answer to the following two natural questions\smallskip

\noindent$\bullet$~Are Chevalley group schemes examples of varieties that can be
defined over $\F_1$?\smallskip

\noindent$\bullet$~Is the notion of variety over $\F_1$ as in \opcit compatible
with the geometric viewpoint developed by Tits?\smallskip

The first question was formulated in \cite{Soule}. In this article we  show
that Chevalley group schemes can be defined\footnote{as varieties but not as groups} over the quadratic extension $\F_{1^2}$
(\cf~Theorem~\ref{mainthm}). The second question originates naturally by
working with the simplest example of a projective variety, namely the
projective spaces $\P^d$.

At  first sight, a very serious problem emerges in \cite{Soule}, since the
definition of the set $\P^d(\F_1)$  does {\it not} appear to be naturally
linked with the notion of a geometry over $\F_1$ as in \cite{Tits}. In fact, in
\S~6 of \cite{Soule} the cardinality of the set $\P^d(\F_{1^n})$  is shown to
be $N(2n+1)$, where $N(x)$ is a polynomial function whose values at prime
powers $q=p^r$ are given by the classical formula $N(q)= q^d+\ldots + q + 1$
giving the cardinality of $\P^d(\F_{q})$. When $n=1$, one obtains the integer
$|\P^d(\F_{1})|=N(3)=\frac{3^{d+1}-1}{2}$ which is incompatible with (and much
larger than) the number $d+1$ of points of the set $\mathcal P_d$ on which Tits
defines his notion of projective geometry of dimension $d$ over $\F_1$
(\cf~\cite{Tits}, \S~13 p. 285).

\smallskip
After clarifying a few statements taken from \cite{Soule} on the notion of
variety over $\F_1$  and on the meaning of a natural transformation of functors
(\cf~\S~\ref{varoverf1}), we show how to resolve the aforementioned problem by
introducing a suitable refinement of the notion of affine algebraic variety over
$\F_1$. The main idea is to replace the category of sets (in which  the
covariant functor $\underline X$ of \opcit takes values)  by the category of
$\Z_{\ge 0}$-{\it graded sets}. In \S~\ref{gadgets}
(\cf~Definition~\ref{gradedgadget}), we explain how to refine the covariant
functor $\underline X$ into a graded functor $\underline X=\coprod_k \underline
X^{(k)}$ defined by a disjoint union of homogeneous components which
correspond, at the intuitive level, to the terms of the Taylor expansion, at
$q=1$, of the counting function $N(q)$. The condition that $N(q)$ is a
polynomial is very restrictive (\eg it fails in general for elliptic curves) and was
required in \cite{Soule} (\S 6, Condition Z) to define the zeta function\footnote{The definition of the zeta function given there is ``upside down" with respect to Manin's, and should be replaced by its inverse to get \eg $\zeta_{\P_1}(s)=1/(s(s-1))$ instead of $s(s-1)$.   } .  In
\S~\ref{elexamples} we check that in the case of a projective space, the set
$\underline \P^d(\F_{1^n})$ coincides {\em in degree zero} with the $d+1$
points of the set $\mathcal P_d$, and this result shows the sought for
agreement with the theory of  Tits.

 Our new definition of an affine variety over $\F_1$ is described by the following data

 (a)~A  covariant functor from the category of finite abelian groups
  to the category of graded sets
 $$\underline X=\coprod_{k\ge 0} \underline X^{(k)}: \cF_{ab}\to \mathcal Sets.$$

 (b)~An affine  variety $X_\C$ over $\C$.

 (c)~A natural transformation $e_X$ connecting $\underline X$
 to the functor$$\cF_{ab}\to \mathcal Sets,\quad\abgr\mapsto \Hom(\Spec\, \C[\abgr],X_\C).$$

These data need to fulfill also a strong condition
(\cf~Definition~\ref{affvarfun}) which determines uniquely a variety over $\Z$.
To a point of $\Spec\, \C[\abgr]$ is associated a character $\chi: \abgr\to
\C^*$ which assigns to a group element $g\in \abgr$ a root of unity $\chi(g)$
in $\C$.  For each such character, the map $e_X$ provides a concrete
interpretation of the elements of $\underline X(\abgr)$ as points of $X_\C$.

\smallskip
In \S~\ref{Chevgadg} we test these ideas with Chevalley groups $G$ and show
that they can be defined over $\F_{1^2}$. Let $\mathfrak G$ be the algebraic
group scheme over $\Z$ associated by Chevalley (\cite{Chevalley1},
\cite{Demazure}) to a root system $\{L,\Phi,n_r\}$ of $G$ (\cf~\S~\ref{rootsys}
and \cite{Tits1}, \S~4.1). In \S~\ref{theproof} we prove that $G$ can be
defined over $\F_{1^2}$ in the above sense. For the proof, one needs to verify
that the following conditions are satisfied\smallskip

\noindent$-$~The functor $\underline G$ (to graded sets) contains enough points
  so that, together with $G_\C$, it characterizes $\mathfrak G$.\smallskip

\noindent$-$~The cardinality of $\underline G(\F_{1^n})$ is given by a
  polynomial $P(n)$ whose value, for $q$ a prime power and $n=q-1$, coincides with
 the cardinality of $\mathfrak G(\F_q)$.\smallskip

\noindent$-$~The  terms of lowest degree in $\underline G$ have degree equal to
the
  rank of $G$ and determine the group extension of the Weyl group of $G$ by $\Hom(L,\abgr)$, as
  constructed by Tits in \cite{Tits1}.\smallskip

The first condition ensures the compatibility with Soul\'e's original notion of
variety over $\F_1$. The second statement sets a connection with the theory of
zeta functions as in  \cite{Manin}. Finally, the third condition guarantees a
 link with the constructions of Tits. In fact, in \cite{Tits} it was originally promoted the idea that the Weyl group of a Chevalley group $G$ should be interpreted as
the points of $G$ which are rational over $\F_1$. For $G=\GL_{d}$, it was then
shown in \cite{Kapranov} that the points of $G$ over $\F_{1^n}$ are described
by the wreath product of  the group of permutations of $d$ letters with
$\mu_n^d$. It is important to notice that in our theory these groups are recast
as  the terms of lowest degree of $\underline G$. The terms of higher degree
are more subtle to describe; to construct them we  make use of the detailed
theory of Chevalley as in ~\cite{Chevalley}, \cite{Chevalley1}.

If $G$ is the  Chevalley group associated to a root system, the cardinality of the set $\mathfrak G(\F_q)$ (\ie the number of
points of $\mathfrak G$ which are rational\footnote{There are in general more
points in $\mathfrak G(\F_q)$ than in the Chevalley group  $G_{\F_q}$ which
is the commutator subgroup of $\mathfrak G(\F_q)$.} over $\F_q$), is given by
the formula
\begin{equation}\label{chevgroup}
|\mathfrak G(\F_q)| = (q-1)^\ell q^N \sum_{w\in W}q^{N(w)},
\end{equation}
where $\ell$ denotes the rank of $G$, $N$ is the number of positive roots,
 $W$ is the Weyl group and $N(w)$ is
the number of positive roots $r$, such that $w(r)<0$. The above formula
\eqref{chevgroup} corresponds to a decomposition of $\mathfrak G(\F_q)$ as a
disjoint union (over the Weyl group $W$) of products of the form
\begin{equation}\label{gfq}
    \mathfrak G(\F_q)=\displaystyle{\coprod_{w\in W}}\,\bG_m(\F_q)^\ell\,\A^N(\F_q)\, \A^{N(w)}(\F_q).
\end{equation}
This equality suggests the definition of the functor $\underline G$  by means
of a sum  of products of powers of the graded functors $\underline \bG_m$ and
$\underline \A$ (\cf~\S~\ref{elexamples}).

The most technical part of our construction is the definition of the natural
transformation $e_G$ as in  (c), which involves the introduction of a lifting
procedure from the Weyl group $W$ to the complex group $G_\C=\mathfrak G(\C)$. The solution to
this problem is in fact already contained in \cite{Tits1}. Indeed, Tits
introduced in that paper  a functor $\cN_{\abgr,\epsilon}(L,\Phi)$ from pairs
$(\abgr,\epsilon)$ of an abelian group and an element of square one in $D$, to
group extensions of the form
$$1\to
\Hom(L,D)\to\cN_{\abgr,\epsilon}(L,\Phi)\to W\to 1.$$ The
 definition of the graded functor $\underline G$ and the natural
  transformation $e_G$ then follow by applying the original
  construction of Tits together with the Bruhat decomposition
  of $G$ and working with a specific parametrization of its cells.

  Incidentally, we find rather remarkable that the
image of this lift of the Weyl group in $G_\C$ consists only of finite products
of elements in the maximal torus with elements  of the form $x_r(\mu)$, where
$\mu$ is a root of unity in $\C$ and where the $x_r(t)$ generate (over any
field $k$) unipotent one-parameter subgroups  associated to the roots
$r$. The fact that $\underline G$ contains enough points so that, together with
$G_\C$, it characterizes $\mathfrak G$, follows from an important result of
Chevalley \cite{Chevalley1}, by working only with the points in the big cell of
$G$.

The original definition of varieties over $\F_1$ in \cite{Soule}, as well as the variant used here, is based on a covariant construction of enough points with cyclotomic coordinates, however a precise control on the size of this set is also needed.
In \S \ref{finrem}, we show that the above examples of varieties $X$ over $\F_1$ (more precisely $\F_{1^2}$) fulfill much stronger properties than those originally required. That is:\vspace{.05in}

$\bullet$~The construction of the functor
$\underline X$ extends to the category
 of pairs $(\abgr,\epsilon)$ of an arbitrary abelian semi-group\footnote{with a unit and a zero element} and an element of square one.\vspace{.05in}
 
$\bullet$~The definition   of the natural
 transformation $e_X$ extends to arbitrary commutative rings $A$ and determines a
 map $  e_{X,A}\,:\; \underline X(A_s)\to X(A)\,,
$ where $A_s$ is the abelian semi-group given by the multiplication in $A$.\vspace{.05in}

$\bullet$~The natural transformation $ e_{X,A}$ is a bijection when $A=K$ is a field.\vspace{.05in}

In the final part of the paper (\cf~\S~\ref{finrem}),  we explain how this enriched construction  yields  a new notion of schemes over $\F_1$ which reconciles Soul\'e's original viewpoint with the approach taken up by Deitmar in \cite{deit} (following Kurokawa, Ochiai, and Wakayama \cite{KOW}) and the log-geometry of monoids of K.~Kato \cite{Kato}, with also the extra advantage of not being limited to the case of toric varieties. In the Noetherian case, the local representability of the functor $\underline X$ implies that its restriction to finite abelian groups is automatically a functor to graded sets, thus clarifying the role of the grading in our construction.  This new notion of schemes over $\F_1$ supplies also a conceptual reason for the
 equality of the number of points of  $X(\F_q)$
 and  the cardinality  of the set $\underline X(\abgr)$ ($\abgr=\F_q^*$) and for the polynomial nature of the counting function (\cf \cite{ak1}).

\medskip

{\bf Acknowledgment.} The authors are partially supported by the NSF grant
DMS-FRG-0652164. We are grateful to P. Cartier and C. Soul\'e for useful comments.\medskip

\section{On the notion of variety over $\F_1$}\label{varoverf1}

In this section we review the notion of variety over $\F_1$ as in \cite{Soule}
and develop a refinement of this concept that will be applied to the case of
Chevalley group schemes in \S\ref{Chevgadg}, to show that these varieties can
be defined over (an extension of) $\F_1$, thus establishing a link with Tits'
geometries.

\subsection{Extension of scalars}\label{extofscalars}

In this paragraph we shall use the same notation as in \cite{Soule}. Let $k$ be a
field and $\Omega$ be a commutative $k$-algebra. One considers the functor of
extension of scalars:
\begin{equation}\label{extsc}
 \cdot\otimes_k\Omega: \cA_k \longrightarrow \cA_\Omega,\qquad  R\mapsto R_\Omega=R\otimes_k\Omega
\end{equation}
from the category $\cA_k$ of unital commutative $k$-algebras to the
corresponding category $\cA_\Omega$. This functor extends to the category of
schemes over $k$ and we use the same notation to denote it. If $X$ is a scheme
(of finite type) over $k$, one lets $X_\Omega=X\otimes_k\Omega$ the
corresponding scheme over $\Omega$.  If $X=\Sp(R)$ is affine and corresponds to
the $k$-algebra $R$, then $X_\Omega$ is also affine and corresponds to the $
\Omega$-algebra $R_\Omega$. The natural homomorphism of algebras $R\to
R_\Omega$ corresponds at the level of schemes to a surjective morphism
\begin{equation}\label{schfunc}
    X_\Omega=\Sp(R_\Omega)\to X=\Sp(R)\,.
\end{equation}
Let $\mathcal Set$ be the category of sets. Then we view a scheme $X$ over $k$
as a covariant functor
\begin{equation}\label{funsets}
 \underline X\; :\; \cA_k \longrightarrow \mathcal Sets,\quad R\mapsto \underline X(R).
\end{equation}
For affine schemes $X=\Sp(A)$, one has $\underline X(R)= \Hom(A,R)$.  Note that
the functor $X\to \underline X$ on schemes is covariant.\vspace{.05in}

In \cite{Soule} (\cf~Proposition 1), one makes use of the following statement

\begin{prop} \label{basicprop} (i) There exists a natural transformation of functors
\begin{equation}\label{nattrans}
  i\;:\; \underline X\to \underline X_\Omega\,,\ \ \   X(R)\subset X_\Omega(R_\Omega)
\end{equation}
(ii) For any scheme $S$ over $\Omega$ and any natural transformation
\begin{equation}\label{nattrans1}
  \varphi\;:\; \underline X\to \underline S\,,
\end{equation}
there exists a unique morphism  $\varphi_\Omega$ (over $\Omega$) from
$X_\Omega$ to $S$ such that $\varphi=\underline \varphi_\Omega\circ i$.
\end{prop}

Notice that \eqref{nattrans} seems to imply that by means of the covariance
property of the functor $X\mapsto\underline{X}$, one should obtain a natural
morphism of schemes $X\to X_\Omega$, and this is in evident contradiction with
\eqref{schfunc}. This apparent inconsistency is due to an abuse of notation and
it is easily fixed as follows. The functors $\underline X_\Omega$ and
$\underline S$ which in \opcit are defined as functors from $\cA_k$ to
$\mathcal Sets$ (\cf~equation above Proposition 1),  should instead be properly
introduced as
  functors from the category $\cA_\Omega$ to $\mathcal Sets$. The ``hidden'' operation in \cite{Soule} is the composition with the  functor of extension of scalars
  \begin{equation}\label{beta}
    \beta: \cA_k \longrightarrow \cA_\Omega,\quad \beta(R)=R_\Omega
  \end{equation}
Up-to replacing in the above proposition $\underline X_\Omega$ by $\underline
X_\Omega\circ \beta$ and $\underline S$ by $\underline S\circ \beta$, one then
obtains the following correct statement

\begin{prop} \label{basicpropbis} (i) There exists a natural transformation of functors
\begin{equation}\label{nattransbis}
  i\;:\; \underline X\to \underline X_\Omega\circ \beta\,,\ \ \   \underline X(R)\subset \underline X_\Omega(R_\Omega)
\end{equation}
(ii) For any scheme $S$ over $\Omega$ and any natural transformation
\begin{equation}\label{nattrans1bis}
  \varphi\;:\; \underline X\to \underline S\circ \beta\,,
\end{equation}
there exists a unique morphism  $\varphi_\Omega$ (over $\Omega$) from
$X_\Omega$ to $S$ such that $\varphi=\underline \varphi_\Omega\circ i$.
\end{prop}

The proof is a simple translation of the one given in \cite{Soule}, with only a
better use of notation.

\proof We first consider the case of an affine scheme $X=\Sp(A)$. The proof in that case applies to any functor
$\underline S$   from $\cA_\Omega$ to $\mathcal Sets$. The functor $X_\Omega$ is represented by $A_\Omega$ and the inclusion $X(R)\subset X_\Omega(R_\Omega)$ is simply described by the
inclusion
$$
i\; :\; \Hom_k(A,R)\subset \Hom_\Omega(A\otimes_k\Omega,R\otimes_k\Omega)\,, \
\ \ f\mapsto i(f)=f\otimes_k id_\Omega=\beta(f)\,.
$$

Next, since $\underline X$ is represented by $A$, Yoneda's lemma shows that the natural transformation  $\varphi$ of
\eqref{nattrans1bis} is characterized by
\begin{equation}\label{idmorp}
h=\varphi(A)(id_A)\in \underline S(A_\Omega)
\end{equation}
(in the displayed formula just after Proposition 1 in \cite{Soule} there is a
typo: the term $X_\Omega(A_\Omega)$ should be replaced by $\underline
S(A_\Omega)$). Similarly, the functor $X_\Omega$ is represented by $A_\Omega$, and a morphism of
functors $\psi$ from $X_\Omega$ to $\underline S$ is uniquely determined by
 an element of $\underline S(A_\Omega)$. Thus $h$ determines  a unique morphism $\psi=\varphi_\Omega$ from $X_\Omega$ to $\underline S$ such that
 $$
 \psi(A_\Omega)(id_{A_\Omega})=h\,.
 $$
 \smallskip
The equality $\varphi=\underline \varphi_\Omega\circ i$ follows again from Yoneda's lemma since
$\underline X$ is represented by $A$ and one just needs to check the equality on  $id_A$ and it follows from  $i(id_A)=id_{A_\Omega}$.
Similarly, to prove the uniqueness, since the functor $X_\Omega$ is represented by $A_\Omega$, it is enough to show that the equality $\varphi =\psi\circ i$ uniquely determines $\psi_{A_\Omega}(id_{A_\Omega})$.

The extension to schemes which are no longer affine follows as in
\cite{Soule}, but the proof given there is unprecise since it is not true in general that the inverse image of an open affine subscheme of a scheme $X$ by a morphism
$\Spec A\to X$ is affine.  As a functor from $\cA_\Omega$ to $\mathcal Sets$, $\underline X_\Omega$ is the composition $\underline X\circ \alpha$ of $\underline X$ with the restriction of scalars $\alpha$ from $\Omega$ to $k$ which is the right adjoint of $\beta$. The inclusion $i$ comes from the canonical morphism $R\to \alpha(\beta(R))$ for any object $R$ of $\cA_k$. We leave it to the reader to clarify the proof and show the result using on the functor $\underline S$ the only assumption that it is {\em local} in the sense of \cite{demgab} Definition 3.11.
\endproof

\medskip

\subsection{Gadgets}\label{gadgets} We keep the same notation as in the previous paragraph.
Let us first recall the definition of ``truc'' given in \cite{Soule}.
Let $\cR$ be the full subcategory of the category of commutative rings whose objects are the group rings $\Z[H]$ of finite abelian groups (\cf \cite{Soule} 3. D\'efinitions, Remarques).

\begin{defn} \label{truc}
A truc over $\F_1$ consists of the following data\vspace{.05in}

-~a pair $X=(\underline X, \cA_X)$ of a covariant functor $\underline X:
\cR\longrightarrow \mathcal Sets$ and a $\C$-algebra $\cA_X$

-~for each object $R$ of $\cR$ and each homomorphism $\sigma:R\to \C$, an
evaluation morphism ($\C$-algebra homomorphism)
   $$e_{x,\sigma}:\cA_X\to \C\qqq x\in \underline X(R),$$
   which satisfies the functorial compatibility
    $e_{f(y),\sigma} = e_{y,\sigma\circ f}$, $\forall$ $f:R' \to R$
     morphism in $\cR$ and $\forall~y\in \underline X(R')$.
\end{defn}

We shall reformulate slightly the above definition with the goal
to:\vspace{.05in}

$\bullet$~treat the archimedean place simultaneously with $\Sp \,\Z$

$\bullet$~replace the category $\cR$ by the category $\cF_{ab}$  of   finite
abelian groups\footnote{This replacement was suggested by C. Soul\'e}

$\bullet$~put in evidence the role of the functor $\beta$.\vspace{.05in}

Thus, we replace $\cR$ by the category $\cF_{ab}$  of   finite abelian groups.
There is a natural functor of extension of scalars from $\F_1$ to $\Z$ which is
given by
\begin{equation}\label{extofscal}
\beta\;:\; \cF_{ab}\longrightarrow \cR,\quad \beta(\abgr)=
\abgr\otimes_{\F_1}\Z:=\Z[\abgr]
\end{equation}
and associates to an abelian group its convolution algebra over $\Z$. Let us
understand the evaluation morphism as a natural transformation. We introduce
the functor
\begin{equation}\label{archfunct}
 \Sp_\infty(\cA_X)\,: \cR\longrightarrow\mathcal Sets,\quad \ R\mapsto   \Hom(\cA_X, R\otimes_\Z\C)
\end{equation}
and compose it with the functor $\beta: \cF_{ab}\to\cR$.

\begin{lem}\label{naturaltr}
The evaluation morphism $e$ as in Definition~\ref{truc}, determines a natural
transformation of (covariant) functors
\[
e: \underline X \to\Sp_\infty(\cA_X)\circ \beta.
\]
\end{lem}

\proof For each object $\abgr$ of $\cF_{ab}$ the evaluation map $e_{x,\sigma}$
can be viewed as a map of sets $\underline X(\abgr) \to \Hom(\cA_X,
R\otimes_\Z\C)$, where $R= \beta(\abgr)$. Indeed, for $x\in \underline
X(\abgr)$, we get a map $e_x$ from characters $\sigma$ of $R$ to characters of
$\cA_X$ which determines a morphism from $\cA_X$ to $R\otimes_\Z\C$.
\endproof
We now reformulate Definition \ref{truc} as follows

\begin{defn} \label{gadget}
A gadget over $\F_1$ is a triple $X=(\underline X,X_\C,e_X)$ consisting
of\vspace{.05in}

(a)~a  covariant functor $\underline X: \cF_{ab}\to \mathcal Sets$ to the
category of  sets

(b)~a  variety $X_\C$ over $\C$

(c)~a natural transformation $e_X: \underline X \to \Hom(\Spec\, \C[-],X_\C)$
from the functor $\underline X$
 to the functor
 \begin{equation}\label{abgrgad}
 \Hom(\Spec\, \C[-],X_\C),\quad \abgr\mapsto \Hom(\Spec\, \C[\abgr],X_\C)\,.
 \end{equation}
\end{defn}\medskip

\subsubsection{Example}\label{exgadget}
An affine variety $V$ over ${\mathbb Z}$ defines a gadget $X=\cG(V)$ over
${\mathbb F}_1$ by letting $X_\C=V_\C = V\otimes_\Z\C$.  $\underline X (\abgr)=\Hom(\cO,\Z[\abgr])$
is the set of points of $V$ in the convolution algebra  $\Z[\abgr]$ with the
natural transformation to $\Hom(\Spec\, \C[\abgr],V_\C)=\Hom(\cO_\C, \C[\abgr])$ obtained by applying the functor $\otimes_\Z\C$.

\begin{defn} \label{gradedgadget} A gadget $X$ over $\F_1$ is said to be graded when
$$\underline X=\coprod_{k\ge 0} \underline X^{(k)}: \cF_{ab}\to \mathcal Sets$$  takes values in the category of $\Z_{\ge 0}$-graded sets. It is finite when
 the set
$\underline X(\abgr)$ is  finite  $\forall~\abgr\in \cF_{ab}$.
 \end{defn}\medskip

\subsection{Varieties over $\F_1$}
The notion of morphism of gadgets $\phi: X \to Y$ is essentially that of a
natural transformation. More precisely, $\phi$ is determined by a pair
$\phi=(\underline \phi, \phi_\C)$, with
\begin{equation}\label{fungadg}
    \underline \phi: \underline X\to \underline Y\,, \qquad \phi_\C: X_\C\to Y_\C.
\end{equation}
$\underline \phi$ is a natural transformation of functors and $\phi_\C$ a
morphism of varieties over $\C$. A {\em required} compatibility with the evaluation
maps gives rise to a commutative diagram
\begin{gather}
\label{functmap}
 \,\hspace{50pt}\raisetag{-47pt} \xymatrix@C=25pt@R=25pt{
 \underline X(\abgr)\ar[d]_{e_X(\abgr)}\ar[r]^-{\underline \phi(\abgr)} &
  \underline Y(\abgr)\ar[d]^{e_Y(\abgr)}& \\
\Hom(\Spec\, \C[\abgr],X_\C) \ar[r]^-{\phi_\C}  & \Hom(\Spec\, \C[\abgr],Y_\C).\\
}
\end{gather}

\medskip

As in \cite{Soule}, we introduce the following  notion of  {\em immersion} of
gadgets

\begin{defn}\label{injgadg} A morphism of gadgets $\phi: X \to Y$ is said to be an immersion if $\phi_\C$ is an embedding and for any object $\abgr$ of $\cF_{ab}$, the
map $ \underline \phi: \underline X(\abgr)\to \underline Y(\abgr) $ is
injective.
\end{defn}

We can now re-state the key definition  of an {\em affine} variety $X$ over
$\F_1$. In the formulation given in \cite{Soule}, it postulates the existence
of a variety (of finite type) over $\Z$ which plays the role of the scheme
$X\otimes_{\F_1} \Z$ and fulfills the universal property of Proposition
\ref{basicpropbis}.

\begin{defn}\label{affvarfun}
An  {\em affine} variety $X$ over $\F_1$ is a finite, graded gadget $X$ such
that there exists an affine variety $X_{\mathbb Z}$ over ${\mathbb Z}$ and an
immersion $i : X \rightarrow \cG(X_{\mathbb Z})$ of gadgets satisfying the
following universal property: for any affine variety $V$ over ${\mathbb Z}$ and
any morphism of gadgets $ \varphi : X \rightarrow \cG(V)\, , $ there exists a
unique algebraic morphism
$$
\varphi_{\mathbb Z} : X_{\mathbb Z} \rightarrow\textbf{} V
$$
of affine varieties such that $\varphi = \cG(\varphi_{\mathbb Z}) \circ i$.
\end{defn}\medskip

\subsection{Varieties over $\F_{1^n}$} This small variant is obtained
(following \cite{Soule} \S 3.8.2) by replacing the category $\cF_{ab}$ of
finite abelian groups by the finer one $\cF_{ab}^{(n)}$ whose objects are
 pairs $(\abgr,\epsilon)$, where $\abgr$ is a finite abelian group and
$\epsilon\in \abgr$  is of order exactly $n$. A morphism in $\cF_{ab}^{(n)}$
is a homomorphism of abelian groups which sends $\epsilon\in \abgr$ to
$\epsilon'\in \abgr'$. Let $R_n=\Z[T]/(T^n-1)$ then the whole discussion takes place over $\Spec R_n$. We shall only use in this paper the case $n=2$. In that case, the two homomorphisms
$\rho_\pm:R_2\to \Z$ given by $\rho_\pm(T)=\pm 1$ show that $\Spec R_2$ is the union of two copies $(\Spec \Z)_\pm$ of $\Spec \Z$ which cross at the prime $2$. We shall concentrate on the non-trivial copy
$(\Spec \Z)_-\subset \Spec R_2$. In Definition~\ref{gadget} one replaces $\cF_{ab}$ by
$\cF_{ab}^{(2)}$ and  one substitutes everywhere the group ring
$\beta[\abgr]$ by the reduced group ring which is the tensor product of rings
\begin{equation}\label{redgring}
    \beta[\abgr,\epsilon]=\Z[D]\otimes_{\Z[\Z/2\Z]}\Z\,,
\end{equation}
in which $\epsilon=-1$.
Thus the characters $\chi$ of the algebra $\C[\abgr,\epsilon]$ are the
characters of $\C[\abgr]$ such that $\chi(\epsilon)=-1$. They still separate the elements of $D$:
\begin{equation}\label{spchar}
    \forall g_1\neq g_2\in \abgr\,, \ \ \exists \chi\in \Spec\C[\abgr,\epsilon]\,, \ \chi(g_1)\neq \chi(g_2)
\end{equation}
An affine variety $V$ over ${\mathbb Z}$ defines a gadget $X=\cG(V)$ over
${\mathbb F}_{1^2}$ by letting $X_\C=V_\C = V\otimes_\Z\C$ sit over the non-trivial copy
$(\Spec \Z)_-\subset \Spec R_2$. The functor  $\underline X (\abgr)=\Hom(\cO,\beta[\abgr,\epsilon])$
is the set of points of $V$ in the   algebra  $\beta[\abgr,\epsilon])$ with the
natural transformation to $\Hom(\Spec\, \C[\abgr,\epsilon],V_\C)=\Hom(\cO_\C, \C[\abgr,\epsilon])$ obtained by applying the functor $\otimes_\Z\C$.

\medskip

\section{Elementary examples}\label{elexamples}

In this section we apply Definition~\ref{affvarfun} by working out the explicit
description of the graded functor $\underline X$ in several elementary examples
of algebraic varieties over $\F_1$. The main new feature, with respect to
\cite{Soule}, is the introduction of a grading. At the intuitive level, the
underlying principle in the definition of the graded  functor $\underline
X=\coprod_{k\ge 0} \underline X^{(k)}$ is that of considering the Taylor
expansion, at $q=1$, of the function $N(q)$ counting the number of points of
the scheme $X$ over the finite field $\F_q$. The term of degree $k$ (\ie
$a_k(q-1)^k$) in the expansion should agree with the cardinality of the set
$\underline X^{(k)}(\abgr)$, for $q-1=|\abgr|$, $D\in\text{obj}(\cF_{ab})$.

The requirement that the function $N(q)$ counting the number of points of
the scheme $X$ over the finite field $\F_q$ is a polynomial in $q$ is imposed in
\cite{Soule} in order to deal with the zeta function and is very restrictive. It fails for instance in general for elliptic curves but it holds for instance for Chevalley group schemes.
We shall first deal with
 a few concrete examples of simple geometric spaces for which $N(q)$ is easily computable. These are\vspace{.05in}

$(1)$~$\bG_m$, $N(q)=q-1$.

$(2)$~The affine line $\A^1$, $N(q)=q$.

$(3)$~The projective space $\P^d$, $N(q)=1+q+\ldots +q^d$.\vspace{.05in}

In the following we shall consider each of these cases in details.\medskip

\subsection{$\Sp D$} Let $D$ be a finite abelian group. We let $\Sp D$ be the
gadget given by $\underline \Sp D(D')=\Hom(D,D')$ and $(\Sp D)_\C=\C[D]$ with the obvious natural transformation. One checks that
it defines a variety over $\F_1$. It is graded by the grading concentrated in degree $0$.

\subsection{The multiplicative group $\bG_m$}
 For the multiplicative group $\bG_m$, the counting function is $N(q)=q-1$:
its Taylor expansion at $q=1$ has just one term in degree $1$. We define the
functor $\underline \bG_m$ from abelian groups to $\Z_{\ge 0}$-graded sets
accordingly \ie
\begin{equation}\label{bgmbis}
\underline\bG_m: \cF_{ab} \longrightarrow \mathcal Sets,\quad
\underline\bG_m(\abgr)^{(k)}=
    \begin{cases} \emptyset&\text{if}~k\in\Z_{\ge 0}\setminus\{1\}\\
\abgr &\text{if}~k=1.\end{cases}
\end{equation}
In particular one sets
\begin{equation}\label{bgm}
    \underline\bG_m(\F_{1^n})^{(k)}=
    \begin{cases} \emptyset&\text{if}~k\in\Z_{\ge 0}\setminus\{1\}\\
\Z/n\Z &\text{if}~k=1.\end{cases}
\end{equation}
Except for the introduction of the grading and for the replacement of the
category of (commutative) rings finite and flat over $\Z$ (as in \cite{Soule})
by that of finite abelian groups, the definition \eqref{bgmbis} is the same as
the corresponding functor in \opcit

We denote by $e_m: \underline \bG_m \to \Hom(\Spec\, \C[-],\bG_m(\C))$ the
natural transformation from the functor $\underline \bG_m$
 to the functor $$\Hom(\Spec\, \C[-],\bG_m(\C)),\quad \abgr\mapsto \Hom(\Spec\, \C[\abgr],\bG_m(\C)),$$ which assigns to a
 character $\chi$ associated to a point of $\Spec\, \C[\abgr]$ the group homomorphism
 \begin{equation}\label{natransfo}
    \abgr\to \C^*\,,  \ e_m(D)(g)=\chi(g)\,.
 \end{equation}

It is now possible to adapt the proof of \cite{Soule} (as in 5.2.2) and show
that this gadget defines a variety over $\F_1$.

\begin{prop}\label{propmul} The gadget $\bG_m=(\underline\bG_m,\bG_m(\C), e_m)$
defines a variety over $\F_1$.
\end{prop}

 \proof  By construction $\bG_m$ is a  finite and graded gadget.
 It is  easy to guess that $\bG_{m,\Z}=\Sp(\Z[T,T^{-1}])$ while the immersion $i$ is given by the injection $D\to\Hom(\Z[T,T^{-1}],\Z[D])$. Let us see that
the condition of Definition \ref{affvarfun} is also fulfilled. Let $V=\Sp(\cO)$
be an affine variety over $\Z$. Let $ \phi : \bG_m \rightarrow \cG(V)\, , $ be
a morphism of (affine) gadgets. This means that we are given a pair
$(\underline \phi, \phi_\C)$,  where $\phi_\C$ can be equivalently  interpreted
by means of the corresponding homomorphism of $\C$-algebras ($\mathcal O_\C =
\mathcal O\otimes_\Z\C$)
$$
\phi_\C^*\;: \; \cO_\C\to \C[T,T^{-1}].
$$
Furthermore, $\underline \phi$ is a morphism of functors (natural
transformation)
$$
\underline \phi(\abgr)\;:\;  \underline\bG_m(\abgr)\to \Hom(\cO,
\beta(\abgr))\,
$$
which fulfills the following compatibility condition (\cf \eqref{functmap}).
For any finite abelian group $\abgr$ the following  diagram is commutative
\begin{gather}
\label{functmap1mul}
 \,\hspace{50pt}\raisetag{-47pt} \xymatrix@C=25pt@R=25pt{
 \underline \bG_m(\abgr)\ar[d]_{e_m(\abgr)}\ar[r]^-{\underline \phi(\abgr)} &
  \Hom(\cO,\beta(\abgr))\ar[d]^{e_{\mathcal G(V)}(D)}& \\
\Hom(\C[T,T^{-1}],\beta(\abgr)_\C) \ar[r]^-{\phi^*}  & \Hom(\cO_\C,\beta(\abgr)_\C).\\
}
\end{gather}

\bigskip

To construct $\psi=\phi_\Z$ let us show that $\phi^*(\cO)\subset \Z[T,T^{-1}]$. Let $h\in \cO$ and
$f=\phi^*(h)$. Then by construction $f\in \C[T,T^{-1}]$. Let $\abgr=\Z/n\Z$ be
the cyclic group of order $n$ with generator $\xi\in \Z/n\Z$; one has $\xi\in
\underline \bG_m(\abgr)$ and $\underline \phi(\abgr)(\xi)\in
\Hom(\cO,\beta(\abgr))$. By evaluating on $h\in \cO$ one gets $$\underline
\phi(\abgr)(\xi)(h)\in \Z[\abgr]\subset \C[\abgr]\,.$$ It follows from the
commutativity of the diagram \eqref{functmap1mul} that this is the same as
evaluating on $f\in \C[T,T^{-1}]$ the homomorphism $e_m(\abgr)(\xi)$ which
coincides with the quotient map
$$
\theta_n:\,\C[T,T^{-1}]\to \C[\Z/n\Z]\,, \ \ T\mapsto \xi\in \Z/n\Z\,.
$$
This means that $\theta_n(f)\in \Z[\Z/n\Z]$, for all $n$. For $f\in
\C[T,T^{-1}]$ one can compute the coefficient of $T^k$ as the limit
$$
b_k=\lim_{n\to \infty} \frac 1n\,\sum_{a=1}^n f(e^{2\pi i\frac an})e^{-2\pi
ik\frac an}.
$$
When $f(x)=x^m$, the sum $\sum_{a=1}^n f(e^{2\pi i\frac an})e^{-2\pi ik\frac
an}$ is either zero or $n$ and the latter case only happens if $m-k$ is a
multiple of $n$. Thus, $\frac 1n\,\sum_{a=1}^n f(e^{2\pi i\frac an})e^{-2\pi
ik\frac an}$, which only depends on $\theta_n(f)$, is a relative integer if
$\theta_n(f)\in \Z[\Z/n\Z]$. It follows that all the $b_k$ are in $\Z$ and
hence that $f\in \Z[T,T^{-1}]$. Thus $\phi^*(\cO)\subset \Z[T,T^{-1}]$ and this uniquely defines $\psi=\phi_\Z\in \Hom(\bG_{m,\Z},V)$. To check that $\phi=\cG(\psi)\circ i$ one uses the
injectivity of  the map
$\Hom_\Z(\cO,\beta(\abgr))\to \Hom_\C(\cO_\C,\beta(\abgr)_\C)$.  \endproof\medskip

\subsection{The affine space $\A^F$}\label{subsectaff}

For the affine line $\A^1$, the number of points of $\A^1(\F_q)$ is $N(q)=q$.
Thus,
 the Taylor expansion of $N(q)$ at $q=1$ is $q=1+(q-1)$ and has two terms in degree $0$ and
$1$. This suggests to refine the definition of the corresponding functor of
\cite{Soule} as follows. We define  $\underline \A^1$ as the graded functor
\begin{equation}\label{aonebis}
\underline\A^1: \cF_{ab}\longrightarrow\mathcal Sets,\quad
\underline\A^1(\abgr)^{(k)}=
    \begin{cases}  \{0\}&\text{if}~k=0\\
\abgr &\text{if}~k=1\\
\emptyset&\text{if}~k\geq 2.\end{cases}
\end{equation}\vspace{.05in}

More generally,  one may introduce for any finite set $F$ the graded functor
\begin{equation}\label{aonef}
    \underline\A^F(\abgr)^{(k)}=
    \coprod_{Y\subset F,\atop|Y|=k} \;\abgr^Y\,.
\end{equation}
which is just the graded product $(\{0\}\cup \abgr)^F$.

 \begin{prop}\label{defnaF} Let $e_F: \underline \A^F \to \Hom(\Spec\, \C[-],\C^F)$ be the natural transformation from the functor $\underline \A^F$
 to the functor $\abgr\mapsto \Hom(\Spec\, \C[\abgr],\C^F)$ which assigns to a
 point in $\Spec\, \C[\abgr]$, \ie to a character $\chi: \C[\abgr]\to\C^*$, the following map
 \begin{equation}\label{natransf}
    \coprod_{Y\subset F} \;\abgr^Y\to \C^F\,, \quad \ e_F(D)((g_j)_{j\in Y})=(\xi_j)_{j\in F},\quad \xi_j=\begin{cases}\chi(g_j)&\text{if $j\in Y$;}\\
    0&\text{if $j\notin Y$.}\end{cases}
 \end{equation}
 Then, the gadget $\A^F=(\underline\A^F,\C^F,e_F)$
defines a variety over $\F_1$.
 \end{prop}

 The proof is identical to that of \cite{Soule}.

   \subsection{Projective space $\P^d$} In \cite{Soule}, after defining the category $\cA$ of affine varieties over $\F_1$, the general case is obtained using contravariant functors from $\cA$ to $\mathcal Sets$, together with a global $\C$-algebra of functions. In the present paper we deal exclusively with the affine case and only briefly mention how the counting of points is affected by the grading in the case of projective space. We adopt the  Definition \ref{gadget} in the non-affine case. As for schemes, it is natural to require the existence of an open covering by affine open sets $U_\alpha$ such that, on each of them, the subfunctor
   $$
   \underline X_\alpha(\abgr)=\{x\in \underline X(\abgr)\,|\, e_X(x)\in \Hom(\Sp(\C[D],U_\alpha)\}
   $$
of the functor $\underline X$ (together with the affine variety $U_\alpha$ over $\C$ and the restriction of the natural transformation $e_X$) defines an affine variety over $\F_1$. One also requires that the $\underline X_\alpha$ cover $\underline X$ \ie that for any $\abgr$ one has $\underline X(\abgr)=\cup  \underline X_\alpha(\abgr)$. One can then rely on Proposition 5 of \S 4.4 of \cite{Soule} to do the patching. In fact we shall obtain in \S \ref{finrem} a general notion of scheme over $\F_1$. We shall  simply explain here, in the case of projective space $\P^d$,
how to implement  the grading.
More precisely, we define $\underline\P^d$ as the following graded functor
\begin{equation}\label{projdimd}
 \underline\P^d: \cF_{ab}\longrightarrow\mathcal Sets,\quad   \underline\P^d(\abgr)^{(k)}=
    \coprod_{Y\subset \{1,2,\ldots, d+1\}\atop |Y|=k+1} \;\abgr^Y/\abgr,\quad k\ge 0
\end{equation}
where the right action of $\abgr$ on $\abgr^Y$ is the diagonal action.  It follows that the
points of lowest degree in $\underline\P^d(\F_{1^n})$ are simply labeled  by
  $\{1,2,\ldots, d+1\}$. Their number is evidently
   \begin{equation}\label{cardproj}
    \#\,\underline\P^d(\F_{1^n})^{(0)}= d+1.
  \end{equation}
  In particular, this shows that $\underline \P^d(\F_{1^n})$  coincides {\em in degree zero} with the $d+1$ points of the set $\mathcal P_d$ on which Tits defines his notion of projective
geometry of dimension $d$ over $\F_1$.
  It is striking that the right hand side of the formula \eqref{cardproj} is
  independent of $n$. This result is also in agreement with the evaluation at $q=1^n$ of the counting function of the set $\P^d(\F_q)$, namely (with the evaluation at $q=1^n$) of the function $N(q)=\sum_{j=0}^d~q^j$.
\medskip

\section{Chevalley group schemes}\label{Chegrsch}

The main result of this section (Theorem~\ref{mainthm}) is the proof that
Chevalley groups give rise naturally to affine varieties over (an extension of)
$\F_1$. To achieve this result we shall need to apply the full theory of
Chevalley groups both in the classical (\ie Lie-theoretical) and algebraic
group theoretical development.\vspace{.05in}

If $\K$ is an algebraically closed field, a Chevalley group $G$ over $\K$ is a
connected, semi-simple, linear algebraic group over $\K$. By definition of a
linear algebraic group over $\K$, $G$ is isomorphic to a closed subgroup of
some $GL_n(\K)$. The coordinate ring of $G$, as affine linear algebraic variety
over $\K$, is then $\K[G] = \K[x_{ij}, d^{-1}]/I$ \ie a quotient ring of
polynomials in $n^2$ variables with determinant $d$ inverted by a prime ideal
$I$.

As an algebraic group over $\K$, $G$ is also endowed with a group structure
respecting the algebraic structure, \ie $G$ is endowed with the following two
morphisms of varieties over $\K$
\begin{equation}\label{grpstr}
\mu: G\times G \to G,~\mu(x,y)=x y;\quad \iota: G \to G,~\iota(x)=x^{-1}
\end{equation}

Notice that by construction $\K[G]$ is a Hopf algebra  whose coproduct encodes
the group structure. \vspace{.05in}

Let $k$ be the prime field of $\K$ and $K$ an intermediate field: $k\subset
K\subset \K$. Then, the group $G$ is said to be {\em defined} over $K$   if the affine variety $G$ and the group structure are
defined over $K$ and, by extension of scalars, also over any field above it
(\cf~\cite{KM}, Chapter 2 (2.1.1), p. 63). In terms of the Hopf algebra
structure, one is given a Hopf algebra $A$ over $K$ such that, as Hopf algebras
$$
\K[G]=A\otimes_K \K\,.
$$

The property for the group $G$ to be {\em split} over $K$ means that some maximal
torus $T\subset G$ is $K$-isomorphic to $\mathbb G_m\times\cdots\times\mathbb
G_m$ ($d$ copies, $d=\dim~T$) (\cf~\cite{Hu}, Chapter XII, \S\S~34.3, 34.4 p.
219-20).

If $G$ is a linear algebraic group over $K$, the set $$G(K) =
\text{Hom}_K(A,K)$$ is a group called the group of $K$-{\it rational points} of
$G$. One has an identification $G(K) = G\cap \A^{n^2}_K$ using $K$-polynomials
to generate the ideal $I$ (\cf~\opcit Chapter XII, \S~34.1, p.
218).\vspace{.05in}

To a semi-simple, connected algebraic group $G$ defined over $K$ and a $K$-split,
maximal torus $T\subset G$,  one associates the group $\text{Hom}(T,\Gm)$
($\Gm=GL_1(K)$): this  is a free abelian group of rank equal to the dimension
of $T$. The group $\R\otimes_\Z \text{Hom}(T,\Gm)$ plays in this context, the
role of the dual $\mathfrak h^*$ of a Cartan Lie-algebra.

 One shows that there exist sub-tori $S\subset T$, $\dim~S = \dim~T-1$,
 whose centralizers $Z(S)$ are of dimension $\dim~T+2$ and such that
 $Z(S)/S$ is isomorphic either to $SL_2(K)$ or $PGL_2(K)$.
 The study of these groups allows one to introduce pairs of
 elements $\pm\alpha\in\text{Hom}(T,\Gm)$ and by varying $S$
 one defines a full set of roots $\Phi$ (\cf~\opcit Chapter XII, \S~34.5 and \cite{Chevalley2}, \S~25.7).

 If $N$ denotes the normalizer of the  torus $T\subset G$, then the
 (finite) group $W=N/T$ acts on $\R\otimes_\Z \text{Hom}(T,\Gm)$ and is called the ($K$-){\it Weyl group}.

The theory of Chevalley groups over a field $K$ has been further extended in
\cite{Chevalley1} and \cite{Demazure}. To every semi-simple, complex Lie group
$G$, and more   generally to an abstract root system, one associates
canonically a
  group scheme $\mathfrak G$ over $\Z$, such that $G$ gets identified with
the group $\mathfrak G(\C)$ of complex points of $\mathfrak G$. We shall return
to this construction in \S~\ref{Cschemes}.

 \vspace{.05in}

In the next paragraphs we shall first review and then apply a construction due
to M. Demazure and J. Tits (\cf~\cite{Demazure}, \cite{Tits1}) which associates to an algebraic reductive group $G$
defined  over  $K$ and a $K$-split maximal torus $T\subset G$, a {\it
canonical} extension of the Weyl group $W$, obtained by considering the groups
of $K$-rational points of $T$ and of its normalizer.  This construction makes
explicit use of a suitable extension of the Weyl group (so called the extended
Weyl group) whose definition is {\it independent} of the field $K$ and is given
only in terms of the root system of $G$.

The notion of extended Weyl group is related to that of extended Coxeter group
$V$ associated to a Coxeter matrix $M$ and a given abstract root system
$\{L,\Phi,n_r\}$. The group $V$ is  a certain extension of  the Coxeter group
of $M$ by a free abelian group  whose rank equals the cardinality of the set of
the reflections associated to the roots.

\medskip

 \subsection{Root systems and Coxeter groups} \label{rootsys} In this paragraph
 we follow \S~4.1 and \S~2.2 of \cite{Tits1} and we shall review several
 fundamental notions associated to the notion of a root system.

 A root system $\{L,\Phi,n_r\}$ is the data given by: \vspace{.05in}

-~a lattice $L$, \ie a free abelian group of finite rank (the group of
weights);

 -~a finite subset $\Phi\subset L$ (the set of roots);

 -~for each $r\in \Phi$, a $\Z$-valued linear form $n_r: L\to \Z$ (the co-root associated to $r$) \vspace{.05in}

 which satisfy the following axioms:\vspace{.05in}

 $(1)$~$L\otimes \Q$ is generated, as $\Q$-vector space, by $\Phi$ and by the intersection of the kernels of the $n_r$;

 $(2)$~$n_r(r)=2$, $\forall~r\in \Phi$;

 $(3)$~the relations $r\in \Phi$, $ar\in \Phi$, $a \in \Q$ imply $a=\pm 1$;

 $(4)$~if $r,s \in \Phi$, then  $r-n_s(r)s\in \Phi$.\vspace{.05in}

For each $r\in\Phi$, the {\it reflection associated to $r$} is the map $s_r: L
\stackrel{\sim}{\to} L$ defined by $s_r(x) = x - n_r(x)\cdot r$. The following equality holds
$s_r = s_{-r}$. \vspace{.05in}

We recall (\cf~\cite{Demazure}, Expos\'e XXI) that the lattice $L$ can be
endowed with a total order which divides the set $\Phi$ of the roots  into two disjoint sets:
positive and negative roots. The set of positive roots $\Phi^+$ is a subset of
$\Phi$ satisfying the conditions:\vspace{.05in}

-~if $r_1, r_2\in \Phi^+$ then $r_1+r_2\in \Phi^+$

-~for each $r\in \Phi$ exactly one of the conditions $r\in \Phi^+$, $-r\in
\Phi^+$ holds.\vspace{.05in}

One  then lets $\Phi^o\subset \Phi^+$ be the set of simple (indecomposable) roots of
$\Phi^+$. It is a
collection $\Phi^o=\{\rho_i, i\in\Pi\}\subset \Phi$ ($\Pi=$ finite set) of
linearly independent roots such that every root can be written as an integral
linear combination of the $\rho_i$, with integer coefficients either all
non-negative or all non-positive. Then, the square matrix $M = (m_{ij})$, $i,j\in\Pi$,
with $2m_{ij}$ equal to the number of roots which are a linear combination of
the $\rho_i$ and $\rho_j$, is a {\it Coxeter matrix} \ie a symmetric square
matrix with diagonal elements equal to $1$ and off diagonal ones positive
integers $\geq 2$.

The {\it Coxeter group} $W=W(M)$ associated to $M$ is defined by a system of
generators $\{r_i, i\in \Pi\}$ (the {\it fundamental reflections}) and
relations
\begin{equation}\label{coxerelat}
    (r_ir_j)^{m_{ij}}=1\qqq i,j\in \Pi.
\end{equation}
The group generated by the reflections  associated to the roots  is canonically
isomorphic to the Coxeter group $W(M)$: the isomorphism is defined by sending
the reflection associated to a simple root to the corresponding generator, \ie
$s_{\rho_i}\mapsto r_i$.

The root system $\{L,\Phi, n_r\}$ is said to be {\it simply connected} if the
co-roots $n_r$ generate the dual lattice $L'=\Hom(L,\Z)$ of $L$. In this case,
the co-roots $n_{\rho_i}$ determine a basis of $L'$.

Given a root system $\{L,\Phi, n_r\}$, there exists a  simply connected root
system $\{\tilde L,\tilde\Phi, n_{\tilde r}\}$  and a homomorphism $\varphi:L
\to \tilde L$ uniquely determined, up-to isomorphism, by the following two
conditions:
\begin{equation}\label{cover}
    \varphi(\Phi)=\tilde\Phi\,, \ \ n_r=n_{\varphi(r)}\circ \varphi\qqq r\in
    \Phi\,.
\end{equation}
The restriction of $\varphi$ to $\Phi$ determines a bijection of $\Phi$ with
$\tilde \Phi$.\vspace{.05in}

The {\it Braid group} $B=B(M)$ associated to a Coxeter matrix $M$ is defined by
generators $\{q_i, i\in \Pi\}$ and the relations
\begin{equation}\label{braidrel}
   \pp(m_{ij};q_i,q_j):= \underbrace{\cdots q_iq_jq_iq_j}_{m_{ij}}=
   \pp(m_{ij};q_j,q_i),\quad \forall~i,j\in \Pi\,.
\end{equation}

The {\it extended Coxeter group} $V=V(M)$ associated to a Coxeter matrix $M$ is
the quotient of $B(M)$ by the commutator subgroup of the kernel $X(M)$ of the
canonical surjective homomorphism $B(M)\to W(M)$. It is defined by generators
and relations as follows. One lets $S$ be the set of elements of $W$ which are
conjugate to one of the $r_i$, \ie the set of reflections (\cf \cite{Tits1},\S
1.2). One considers two sets of generators:
$$\{q_i, i\in \Pi\}, \quad \{g(s), s\in S\}$$ The relations are given by \eqref{braidrel} and the following:
\begin{enumerate}
    \item $q_i^2=g(r_i)$,\quad  $\forall~i\in \Pi$,
  \item $q_i\cdot g(s)\cdot q_i^{-1}=g(r_i(s))$, \quad $\forall~s\in S$, $i\in \Pi$.
  \item $[g(s),g(s')]=1,\quad\forall s,s'\in S$.
\end{enumerate}\vspace{.05in}

Tits shows in Th\'eor\`eme~2.5 of \cite{Tits1} that the subgroup $U=U(M)\subset
V$ generated by the $g(s)$ ($s\in S$) is a free, abelian, normal subgroup of
$V$ that coincides with the kernel of the natural surjective map $f: V \to W$,
$f(q_i) = r_i$ (and $f(g(s))=1$). The group $U=U(M)$ is the quotient of $X(M)$
by its commutator subgroup. \vspace{.05in}

In the following paragraph we shall recall the construction of \cite{Tits1} of
the extended Weyl group using  the extended Coxeter group $V$.\medskip

  \subsection{The group $\cN_{\abgr,\epsilon}(L,\Phi)$}\label{groupN} We keep the same notation
  as in the previous paragraph.  In \S~4.3 of \cite{Tits1} Tits introduces,
  for a given root system $\{L,\Phi,n_r\}$, a functor  $$(D,\epsilon)\to\{N,p,N_s; s\in S\}=\mathcal N$$
  which associates to the pair $(D,\epsilon)$ of an abelian group and an
  element $\epsilon\in D$ with $\epsilon^2=1$, the data
  (\ie an object $\mathcal N$ of a suitable category) given
  by a group $N=\cN_{\abgr,\epsilon}(L,\Phi)$, a surjective
  homomorphism of groups $p: N \to W = W(M)$ and for each reflection $s\in S$,
  a subgroup $N_s\subset N$ satisfying the following conditions:\vspace{.05in}

  (n1)~$\text{Ker}(p)$ is an abelian group;

  (n2)~For $s\in S$, $n\in N$ and $w = p(n)$, $nN_sn^{-1} = N_{w(s)}$;

  (n3)~$p(N_s) = \{1,s\},~\forall~s\in S$.\vspace{.05in}

   A morphism connecting two objects $\mathcal N$ and $\mathcal N'$ is a
   homomorphism $a: N\to N'$ such that $p'\circ a=p$ and $a(N_s)\subset N'_s$ for all $s\in S$ (\cf~\S~3 of \cite{Tits1}).

   Notice, in particular, that the data  $\{V,f,V_s; s\in S\}$, characterizing the  extended Coxeter group, where $V_s \subset V$ is the subgroup generated by $Q_s = \{v\in V, v^2 = g(s)\}$, satisfy (n1)-(n3).\vspace{.05in}

   The object $\mathcal N$ is obtained by the following canonical construction.
   One considers the abelian group $T=\Hom(L,\abgr)$ endowed with the natural (left) action of $W$
   (denoted by $t\mapsto w(t)$, for $w\in W,\, t\in T$)  induced by the corresponding action
   on $L$ (generated by the reflections associated to the roots).
   For each $r\in \Phi$, let $s=s_r\in W$ be the
  reflection associated to the root $r$.
  One lets $T_s$ be the subgroup of $T$ made by  homomorphisms of the form
  \begin{equation}\label{morph}
    L\ni x\mapsto a^{\nu(x)}
  \end{equation}
for some $a\in \abgr$ and where $\nu:L\to \Z$ is a linear form proportional to
$n_r$. Also, one defines (for $s=s_r$)
\begin{equation}\label{morph1}
    h_s(x)= \epsilon^{n_r(x)},\quad x\in L.
  \end{equation}
  (note that replacing $r\to -r$ does not alter the result since
  $\epsilon^2=1$). The formula \eqref{morph1} determines a map $h: S \to T$. Then, the data $\{T, T_s, h_s; s\in S\}$ fulfill the following conditions  $\forall~w\in W,~s\in S$:\vspace{.05in}

$(1)$~$w(T_s)=T_{w(s)}$;

$(2)$~$w(h_s)=h_{w(s)}$;

$(3)$~$h_s\in T_s$;

$(4)$~$s(t)\,\cdot t^{-1}\in T_s$,\quad $\forall~t\in T$;

$(5)$~$s(t)=t^{-1}$,\quad $\forall~t\in T_s$.\vspace{.05in}

One then obtains  the object $\mathcal N$ as follows (\cf Proposition 3.4 of
\cite{Tits1}). One defines the group $N = \cN_{\abgr,\epsilon}(L,\Phi)$ as the
quotient of the semi-direct product group $T\rtimes V$ ($V=$ extended Coxeter group) by the
graph of the homomorphism $U\to T$ ($U = U(M)$) which extends the map
$g(s)\mapsto h_s^{-1}$. One identifies $T$ with its canonical image in $N$ and
for each $s\in S$, one lets $N_s$ be the subgroup of $N$ generated by the
canonical image of $T_s\times Q_s$, where $T_s$ is as above and $Q_s = \{v\in V, v^2 = g(s)\}$. The surjective
group homomorphism $p: N \to W = W(M)$ is induced by $id\times f$ where $f$ is
the canonical group homomorphism $f: V \to W = W(M)$. By construction one has a
morphism connecting the data $\{V,f,V_s; s\in S\}$ to $\{N,p,N_s; s\in
S\}=\mathcal N$. It gives a homomorphism $a: V \to N$ such that, in particular,
    \begin{equation}\label{extcox}
        a(g(s))=h_s\qqq s\in S\,.
    \end{equation}

\medskip

 More precisely, one has the following result (\cf~\S~3.4
and \S~4.3 of \cite{Tits1})

  \begin{prop}\label{titsns} The data $\mathcal N = \{N, p, N_s; s\in S\}$ satisfy the conditions (n1)-(n3).
   Moreover, every map $$\alpha: \{q_i, i\in\Pi\} \to N$$ such that
   $\alpha(q_i)\in N_{r_i}\setminus T_{r_i} = N_{r_i} \cap p^{-1}(r_i)$
  extends to a homomorphism of groups $V\to N$.
  \end{prop}

   \vspace{.05in}

We collect together, for an easy reference, the main properties of the
construction of \cite{Tits1} reviewed in this paragraph.

  \begin{prop} \label{propgrext}  Let $\{L,\Phi,n_r\}$ be a root system. To a pair $(D,\epsilon)$ of
  an abelian group and an element $\epsilon\in D,~\epsilon^2 =1$, corresponds a
  canonical extension  $\cN_{\abgr,\epsilon}(L,\Phi)$  of the Coxeter group $W$
  by $T=\Hom(L,D)$
  \begin{equation}\label{grext}
    1\to T\to N\stackrel{p}{\longrightarrow}W\to 1
  \end{equation}
  and for each reflection $s\in S$  a subgroup $N_s\subset N$. These data satisfy the following properties:
  \begin{itemize}
  \item $nN_sn^{-1} = N_{w(s)}$, for $s\in S$, $n\in N$ and $w = p(n)$.
    \item $p(N_s)=\{1,s\}$, \quad $\forall~s\in S$.
    \item $N_s\cap T=T_s$,\quad  $\forall~s\in S$.
    \item $a^2=h_s\in T_s$,\quad $\forall~a\in p^{-1}(s)=N_s\setminus T_s$.
    \item For each pair $i\neq j$ in $\Pi$, let $m=m_{ij}$ be the order of $r_ir_j\in
    W$, then
    \begin{equation}\label{coxrel}
       \pp(m;a_i,a_j)=\pp(m;a_j,a_i)\qqq a_k\in N_{r_k}\,, \  p(a_k)=r_k\neq 1\,.
    \end{equation}
    \end{itemize}
    The canonical extension $\cN_{\abgr,\epsilon}(L,\Phi)$  of $W$ by $T$   is
functorial in the pair $(\abgr,\epsilon)$, with respect to  morphisms $\ t:D\to
D'$ such that $t(\epsilon)=\epsilon'$.
    \end{prop}

    The meaning of  equation \eqref{coxrel} is the following one. Once a choice of a section $W\supset\Phi^o\ni s\mapsto \alpha(s)\in N_s$ of the
    map $p$ has been made on the set of simple roots $\Phi^o\subset W$, this section admits a natural extension
    to all of $W$ as follows. One writes $w\in W$ as a word  of minimal length
     $w=\rho_1\cdots \rho_k$, in the
    generators $\rho_j\in \Phi^o$. Then \eqref{coxrel} ensures that the corresponding
    product $\alpha(w)=\alpha(\rho_1)\cdots \alpha(\rho_k)\in N$ is independent of the choice of the word
    of minimal length representing $w$ (\cf \cite{Tits1} Proposition 2.1).

\medskip

\subsection{Chevalley Schemes}\label{Cschemes} We keep the notation of \S~\ref{rootsys} and \S~\ref{groupN}.
To a root system $\{L,\Phi, n_r\}$  one associates,
 following \cite{Chevalley1} and \cite{Demazure}, a reductive group scheme $\mathfrak G =\mathfrak G\{L,\Phi,n_r\}$ over $\Z$: the {\it Chevalley scheme}.
 We denote by $\cT$ a maximal torus that is part of a split structure of $\mathfrak G$ and by $\cN$ its normalizer.

To a reflection $s\in S$  correspond naturally the following data: a one
dimensional sub-torus $\cT_s\subset \cT$, a rank one semi-simple subgroup
$\mathfrak G_s\subset \mathfrak G$ containing $\cT_s$ as
 a maximal torus and a point $h_s\in\cT_s$ (belonging to the center of $\mathfrak
 G_s$). One denotes by $\cN_s$ the normalizer of $\cT_s$ in $\mathfrak G_s$.

Let $A$ be a commutative ring with unit and let $A^*$ be its multiplicative
group. We denote by $\mathfrak G(A)$, $\cT(A)$, \resp $\cN(A)$ the groups of
points of $\mathfrak G$, $\cT$, \resp $\cN$ which are rational over $A$. The
quotient $\cN(A)/\cT(A)$ is canonically isomorphic to $W=W(M)$. More precisely,
there exists a unique surjective homomorphism $p_A: \mathcal N(A) \to W(M)$,
whose kernel is $\mathcal T(A)$ so that the data $\{\mathcal N(A), p_A,
\mathcal N_s(A); s\in S\}$ satisfy the conditions (n1)-(n3) of
\S~\ref{groupN}.\vspace{.05in}

The goal of this paragraph is to review a fundamental result of \cite{Tits1}
which describes the data above only in terms of $A$ and the root system
$\{L,\Phi,n_r\}$. To achieve this result one makes use of the following
facts:\vspace{.05in}

-~The group $\cT(A) $ is canonically  isomorphic to  $\Hom_\Z(L,A^*)$ and the
left action of $W(M) \simeq \cN(A)/\cT(A)$ on $\cT(A)$ is induced from the
natural action of $W$ on $L$.

-~If $s=s_r$ is the reflection associated to a root $r\in \Phi$, then
$$
    \cN_{s_r}(A)\cap \cT(A)=\cT_{s_r}(A)=\{\rho\in \Hom(L,A^*)\,|\, \exists \,a\in A^*\,, \
    \rho(x)=a^{\nu(x)}\qqq x\in L\}
 $$
where $\nu: L \to \Z$ is a linear form proportional to $n_r$
(\cf~\S~\ref{groupN}). Taking $a=-1$ and $\nu=n_r$, one gets the element
$h_s(A)\in \cT_{s}(A)$.

-~The normalizer $\cN_s$ of  $\cT_s$ in $\mathfrak G_s$ is such that all
elements of $\cN_{s}(A)$ which are not in  $\cT_{s}(A)$ have a square equal to
$h_s(A)\in \cT_{s}(A)$.
 \vspace{.05in}

We are now ready to state Theorem 4.4 of \cite{Tits1} which plays a key role in
our construction.

\begin{thm} \label{titshelp} The group extension $$1\to \cT(A) \to\cN(A) \stackrel{p}{\longrightarrow} W\to 1$$  is canonically
isomorphic to the group extension $$1\to
\Hom(L,A^*)\to\cN_{A^*,-1}(L,\Phi)\stackrel{p}{\longrightarrow} W\to 1.$$
\end{thm}\medskip

Here $\cN_{A^*,-1}(L,\Phi)$ refers to the functorial construction of
Proposition \ref{propgrext} for the group $D=A^*$ and $\epsilon=-1\in D$. Note
that the case $A=\Z$  corresponds to $D=\{\pm 1\}$, $\epsilon=-1$, and gives
the extension $\cN(\Z)$ of $W$ by $\Hom(L,\{\pm 1\})\simeq (\Z/2\Z)^\ell$. This
particular case contains the essence of the general construction since for any
pair $(\abgr,\epsilon)$ the group $\cN_{\abgr,\epsilon}(L,\Phi)$ is the
amalgamated semi-direct product
$$
\Hom(L,D)\rtimes_{\Hom(L,\{\pm 1\})}\cN(\Z)\,.
$$

\subsection{Bruhat decomposition}\label{Bdec} We keep the notation as in the earlier paragraphs of this section. To each root
$r\in \Phi$ corresponds a root subgroup $\mathfrak X_r\subset \mathfrak G$
defined as the range of an isomorphism $x_r$ from the additive group
$\bG_{a,\Z}$ to its image in $\mathfrak G$ and fulfilling the equation
\begin{equation}\label{rootequ}
    h x_r(\xi)h^{-1}=x_r(r(h)\xi)\qqq~h\in \cT.
\end{equation}
 We recall the following standard notation and well-known
relations:
\begin{equation}\label{lift1}
  n_r(t) =   x_r(t) x_{-r}(-t^{-1})  x_r(t)\qqq t~\in A^*\,, \ \ n_r=n_r(1)
 \end{equation}
 \begin{equation}\label{lift2}
 h_r(t) =   n_r(t)  n_r(-1)\qqq t~\in A^*\,,h_r=h_r(-1)\,, \  \ h_r(t_1)  h_r(t_2) =   h_r(t_1t_2).
 \end{equation}

Let $r,s\in \Phi$ be linearly independent roots and $t,u\in A$. The {\it
commutator} of $x_s (u)$ and $x_r(t) $ is defined as
\begin{equation}\label{comrel}
[x_s (u),x_r (t)] = x_s (u)^{-1}x_r (t)^{-1}x_s (u)x_r (t).
\end{equation}
The following formula, due to Chevalley, expresses the above commutator   as a
product of generators corresponding to roots of the form $ir+js\in \Phi$, for
$i,j >0$

 \begin{lem}\label{canon4} Let $r, s\in \Phi$ be linearly independent roots. Then
 there exist integers $C_{ijrs}\in \Z$ such that
 \begin{equation}\label{commuroots}
    x_s (u)^{-1}x_r (t)x_s (u)x_r (t)^{-1}=\prod_{i,j}\,
    x_{ir+js} (t^iu^j)^{C_{ijrs}}
 \end{equation}
 where the product is applied to pairs $(i,j)$ of strictly positive integers
 such that $ir+js\in\Phi$, and the terms are arranged in order of increasing $i+j$.
 \end{lem}

 The
above formula holds in particular in the case that $r,s$ belong to the subset
$\Phi^+\subset\Phi$ of the positive roots (\cf~\ref{rootsys} and
\cite{Demazure}, Expos\'e XXI), in which case the product is taken over all
positive roots of the form $ir+js$, $i>0$, $j>0$ in increasing order, for the
chosen ordering of the lattice $L$ (\cf \cite{Demazure}, Expos\'e XXII, Lemme
5.5.6 p. 208).\vspace{.1in}

Let $\cU(A)$ be the subgroup of $\mathfrak G(A) $ generated by the elements
$x_r (t)$ for $r\in \Phi^+$, $t\in A$. By construction, the subgroup
$\cU\subset \mathfrak G$ is generated by the root subgroups $\mathfrak X_r$
corresponding to the positive roots. For any $w\in W$, we let $\Phi_w=\{r\in
\Phi^+|w(r)<0\}$ and we denote by $\cU_w$ the subgroup of $\mathfrak G$ generated by the root
subgroups $\mathfrak X_r$ for $r\in  \Phi_w$.

Chevalley proved in \cite{Chevalley}, Th\'eor\`eme 2 the existence of a
canonical form for the elements of the group $\mathfrak G(K) $, when $K$ is a
field.  We recall this result

\begin{thm}\label{chevdec}
Let $K$ be a field. The group $\mathfrak G(K) $ is the disjoint union of the
subsets (cells)
$$C_w=\cU(K)\, \cT(K)\, n_w\,\cU_w(K) $$ where for each $w\in W$, $n_w\in
\cN(K)$ is a chosen coset representative for $w$. The natural map
\begin{equation}\label{chevdec1}
    \varphi_w\,:\,\cU(K)
\times \cT(K) \times \cU_w(K)\to C_w\,, \ \ \varphi_w(x,h,x')=xhn_wx'
\end{equation}
is a bijection for any $w\in W$.
\end{thm}\medskip

We refer to \cite{Demazure}, Expos\'e XXI, Th\'eor\`eme 5.7.4 and Remarque
5.7.5.

\subsection{Chevalley group schemes as gadgets}\label{Chevgadg}

For the definition of the gadget over $\F_1$ associated to a Chevalley group
$G$ and in particular for the construction of the natural transformation
$e_{G}$ (\cf~\S~\ref{gadgets}), one needs to choose a Chevalley basis of the
(complex) Lie algebra of $G$ and a total ordering of the lattice $L$. We keep the
notation as in \S~\ref{Bdec} and in \S~\ref{rootsys} through \S~\ref{Cschemes}.

Chevalley proved (\cf \cite{Chevalley}, Lemma 6 and \cite{Chevalley1}) that, over any commutative
ring $A$, each element of $\cU(A) $ is {\it uniquely} expressible in the form
\[
\prod_{r\in \Phi^+}x_{r} (t_r)\,, \ \ t_r\in A \qqq r~\in \Phi^+
\]
where the product is taken over all positive roots in increasing order. More
precisely one has the following
 \begin{lem}\label{canon} The map
 \begin{equation}\label{canon1}
    t=(t_r)_{r\in \Phi^+}\mapsto \psi(t)= \prod_{r\in \Phi^+} x_r (t_r)
 \end{equation}
establishes a bijection of the free $A$-module   with basis the positive roots
in $\Phi^+$ with  $\cU(A)$.
 \end{lem}

 The proof of this lemma applies without change to give the following
 variant, where we let $\Phi_w=\{r\in \Phi^+|w(r)<0, w\in W\}$
 and we denote by $\cU_w=\displaystyle{\prod_{r\in \Phi_w}}\mathfrak X_r$.

 \begin{lem}\label{canon2} The map
 \begin{equation}\label{canon3}
    t=(t_r)_{r\in \Phi_w}\mapsto \psi_w(t)=  \prod_{r\in \Phi_w} x_r (t_r)
 \end{equation}
establishes a bijection of the free $A$-module with basis $\Phi_w$ with
$\cU_w(A)$.
 \end{lem}
 The key identity in the proof of Lemmas~\ref{canon} and
\ref{canon2} is the commutator relation of Lemma~\ref{canon4}.\vspace{.05in}

We are now ready to apply the theory reviewed in the previous paragraphs to
construct the functor $$\underline{G}: \cF_{ab}^{(2)} \to \mathcal Sets$$ from
the category
 $\cF_{ab}^{(2)}$ of pairs $(\abgr,\epsilon)$ of a finite abelian group and an element of square one, to
the category of graded sets.

\begin{defn}\label{defG} The functor $\underline{G}: \cF_{ab}^{(2)}\to \mathcal Sets$ is defined as the graded product
\begin{equation}\label{unG}
    \underline{G}(\abgr,\epsilon)=\underline\A^{\Phi^+}(\abgr)\times \coprod_{w\in
    W}(p^{-1}(w)\times \underline\A^{\Phi_w}(\abgr))
\end{equation}
where $p$ is the projection
$\cN_{\abgr,\epsilon}(L,\Phi)\stackrel{p}{\longrightarrow} W$ as in
\S\S\ref{groupN} and \ref{Cschemes}. All elements of $p^{-1}(w)$
    have degree equal to the rank of $\mathfrak G$.
\end{defn}

It follows immediately that there are no elements of degree less than the rank
$\ell$ of $\mathfrak G$ and that the set of elements of degree $\ell$ is
canonically identified with $\cN_{\abgr,\epsilon}(L,\Phi)$.\smallskip

We now move to the definition of the natural transformation $$e_{G}:\underline
{G} \to \Hom(\Spec\, \C[-],G_\C),\quad (\abgr,\epsilon)\mapsto \Hom(\Spec\,
\C[\abgr,\epsilon],G_\C).$$ For this part, we make use of the natural
 transformations $e_F$ of \eqref{natransf} for $F=\Phi^+$ and $F=\Phi_w$ and of Theorem~\ref{titshelp}  to
 obtain, for a given character $\chi$, $\chi(\epsilon)=-1$, associated to a point in $\Spec\, \C[\abgr,\epsilon]$, maps
 \begin{equation}\label{emaps1}
    e_{\Phi^+}\,:\,\underline\A^{\Phi^+}(\abgr)\to \C^{\Phi^+}
 \end{equation}
 \begin{equation}\label{emaps2}
    e_{\Phi_w}\,:\,\underline\A^{\Phi_w}(\abgr)\to \C^{\Phi_w}
 \end{equation}
 \begin{equation}\label{emaps3}
    e_{\cN}\,:\,\cN_{\abgr,\epsilon}(L,\Phi)\to \cN(\C).
 \end{equation}
 The last map is compatible with the projection $p$, thus
 maps $p^{-1}(w)$ to $p^{-1}(w)$ for any $w\in W$. We now make use of Lemmas~\ref{canon} and \ref{canon2} to obtain the natural transformation $e_{G}$ defined as follows
 \begin{equation}\label{egnat}
    e_{G}(a,n,b)=\psi(e_{\Phi^+}(a))\,e_{\cN}(n)\, \psi_w(e_{\Phi_w}(b))\in \mathfrak G(\C)=G_\C
 \end{equation}
where $\psi$ and $\psi_w$ are defined as in Lemmas~\ref{canon} and
 \ref{canon2}, for $A=\C$.

\medskip

\subsection{Proof that $G$ determines a  variety over $\F_{1^2}$}\label{theproof}
In this paragraph we shall prove that the gadget $G=(\underline G,G_\C, e_G)$
associated to a Chevalley group $G$ (or equivalently to its root system
$\{L,\Phi,n_r\}$, \cf~\S~\ref{rootsys}) defines a variety  over
$\F_{1^2}$. We keep the earlier notation. We first recall the following
important result of Chevalley (\cf~\cite{Chevalley1}, Proposition~1).

\begin{prop} \label{bigcell} Let $w_0\in W$ be the unique element of the Weyl group
such that $w_0(\Phi^+)=-\Phi^+$ and let $w'_0$ be a lift of $w_0$ in $G_\Z$.
Consider the following morphism, associated to the product in the group,
\begin{equation}\label{maptheta}
  \theta: \cU \times p^{-1}(w_0)\times \cU \to \mathfrak G,\quad  \theta(u,n,v)=unv
\end{equation}
  Then $\theta$ defines an isomorphism of $\cU \times p^{-1}(w_0)\times \cU $ with
  an open affine (dense) subscheme $\Omega $ of $\mathfrak G $, whose global algebra of coordinates is of the form
\begin{equation}\label{affopen}
  \cO_\Omega =  \cO_{\mathfrak G} [d^{-1}]
\end{equation}
where $d\in\cO_{\mathfrak G} $ takes the value $1$ on  $w'_0$.
\end{prop}
 We refer also to proposition
4.1.2 page 172 in \cite{Demazure}  combined with  proposition 4.1.5.

The next theorem shows that the gadget $G=(\underline G,G_\C, e_G)$ over
$\F_{1^2}$ fulfills the condition of Definition \ref{affvarfun}.

\begin{thm} \label{mainthm}
The gadget $G=(\underline G,G_\C, e_G)$ defines a variety over $\F_{1^2}$.
\end{thm}
\proof By construction $G$ is a finite and graded gadget. It is easy to guess
that the sought for scheme $G_\Z$ over $\Z$  is the Chevalley group scheme $\mathfrak
G$ associated to the root system $\{L,\Phi,n_r\}$ (\cf~eg~\cite{Demazure},
Corollary 1.2, Expos\'e XXV). One has by construction an immersion of gadgets
$G\hookrightarrow \cG( \mathfrak G)$. The injectivity of the map of sets follows from two facts,
\begin{itemize}
  \item The characters of $\C[\abgr,\epsilon]$ separate points of $\abgr$ (\cf \eqref{spchar}).
  \item The uniqueness of the Bruhat decomposition in $G_\C=\mathfrak G(\C)$.
\end{itemize}
It remains to be checked that $G$ fulfills
the universal property of Definition \ref{affvarfun}. Let $V=\Sp(\cO(V))$ be an
affine variety of finite type over $\Z$ and $ \phi : G \rightarrow \cG(V)\, , $
be a morphism of gadgets. This means that we are given a pair $(\underline
\phi, \phi_\C)$ where
$$
\phi_\C\;: \; \cO_\C(V)\to \cO_\C(\mathfrak G)
$$
is a homomorphism of  $\C$-algebras, and $\underline \phi$ is a natural
transformation of functors from pairs $(\abgr,\epsilon)$ to sets
$$
\underline \phi(\abgr,\epsilon)\;:\; \underline G(\abgr,\epsilon)\to \Hom(\cO(V), \beta(\abgr,\epsilon))\,
$$
which satisfies the following compatibility condition (\cf \eqref{functmap}):
for any pair $(\abgr,\epsilon)$ the following  diagram commutes

\begin{gather}
\label{functmap1}
 \,\hspace{50pt}\raisetag{-47pt} \xymatrix@C=25pt@R=25pt{
 \underline G(\abgr,\epsilon)\ar[d]_{e_G(\abgr,\epsilon)}\ar[r]^-{\underline \phi(\abgr,\epsilon)} &
  \Hom(\cO(V),\beta(\abgr,\epsilon))\ar[d]^{e_{\mathcal G(V)}(D,\epsilon)=\subset}& \\
\Hom(\cO_\C(\mathfrak G),\C[\abgr,\epsilon]) \ar[r]^-{\phi_\C}  & \Hom(\cO_\C(V),\C[\abgr,\epsilon]).\\
}
\end{gather}

\bigskip

One needs to show   that $\phi_\C(\cO(V))\subset \cO(\mathfrak G)$. Let $h\in
\cO(V)$ and $f=\phi_\C(h)$.  Then by construction $f\in \cO_\C(\mathfrak G)$.
By Proposition~\ref{bigcell},  the intersection $\cO_\C(\mathfrak G)\cap
\cO_\Omega$ coincides with $\cO(\mathfrak G)$ since $ \cO_\Omega =
\cO_{\mathfrak G} [d^{-1}]$ while elements of $\cO_\C(\mathfrak G)$ have a
trivial pole part in  $d^{-1}$. Thus it is enough to show that the restriction
of $f$ to the open affine subscheme $\Omega\subset\mathfrak G$ belongs to
$\cO_\Omega$, to conclude that $f\in \cO(\mathfrak G)$. Let us choose a lift
$w'_0$ of $w_0$ in $\cN_\Z$. In fact we can more precisely choose a lift $w'_0$
of $w_0$ in $\cN_{\Z/2\Z,\epsilon}(L,\Phi)$, where $\Z/2\Z$ is the group of order two generated by $\epsilon$ and then take the image of $w'_0$ under the
map which sends $\epsilon$ to $-1$. We  have $p^{-1}(w_0)=w'_0\cT$. As in
\cite{Chevalley1} (\cf~\S~4, Proposition 1), the algebra $\cO_\Omega$ is the
tensor product of the following three algebras:\vspace{.05in}

-~$\cO(\cU)$ which is the algebra of polynomials over $\Z$ generated by
  the coordinates $t_r$ of Lemma \ref{canon}.

-~$\cO(\cT)=\Z[L]$ the group ring of the abelian group $L$.

-~Another copy of $\cO(\cU)$.\vspace{.05in}

 We consider elements of
$\underline G(\abgr,\epsilon)$ of the form
\begin{equation}\label{eltsbigcell}
    g\in C=\underline\A^{\Phi^+}(\abgr)\times p^{-1}(w_0)\times \underline\A^{\Phi^+}(\abgr)
\end{equation}
and use the choice of $w'_0\in p^{-1}(w_0)\subset \cN_{\abgr,\epsilon}(L,\Phi)$
to identify the cosets $p^{-1}(w_0)=\Hom(L,\abgr)w'_0$. Then, we choose
generators $v_j$, $1\leq j\leq \ell$, of the free abelian group $L$ and use
them to identify $\Hom(L,\abgr)$ with the set of $(y_j)_{j\in
\{1,\ldots,\ell\}}$, $y_j\in \abgr$. Each map $y$ from $Y=\Phi\cup
\{1,\ldots,\ell\}$ to $\abgr$ defines uniquely an element $g(y)\in C$
by\footnote{not all elements of $C$ are of this form}
\begin{equation}\label{eltsbigcell1}
    g(y)=(y_r)_{r\in \Phi^+}\times (y_j)_{j\in \{1,\ldots,\ell\}}\times (y_{-r})_{r\in \Phi^+}
\end{equation}
 Then $g(y)\in \underline G(\abgr,\epsilon)$
and $\underline \phi(\abgr,\epsilon)(g(y))\in \Hom(\cO(V),\beta(\abgr,\epsilon))$ so that by
evaluating on $h\in \cO(V)$ one gets
\begin{equation}\label{map}
\underline \phi(\abgr,\epsilon)(g(y))(h)\in \Z(\abgr,\epsilon)\subset \C[\abgr,\epsilon]\,.
\end{equation}
  By
the commutativity of the diagram \eqref{functmap1}, this is the same as
evaluating on $f\in \cO_\C(\mathfrak G)$ the homomorphism $e_G(\abgr,\epsilon)(g(y))$.
We denote by $k=f|_\Omega$ the restriction of $f$ to $\Omega$, it is given by a
polynomial with complex coefficients
\begin{equation}\label{fresom}
    k=P(t_r,u_j,u_j^{-1})\in \C[t_r,u_j,u_j^{-1}]\,, \ r\in \Phi\,, \ 1\leq j\leq \ell\,.
\end{equation}
 Let $n$ be an integer and
$\abgr=(\Z/n\Z)^Y\times \Z/2\Z=\abgr_1\times \Z/2\Z$ be the group $\abgr_1$ of maps from $Y=\Phi\cup
\{1,\ldots,\ell\}$ to the cyclic group of order $n$, times the cyclic group of
order two with generator $\epsilon$  ($\epsilon^2=1$). We denote by $\xi$ the
generator of $\Z/n\Z$ and for $s\in Y$ we let $\xi_s\in \abgr_1$ have all its
components equal to $0\in \Z/n\Z$ except for the component at $s$ which is
$\xi$. One has a homomorphism of algebras
$$
\theta_n:\,\C[t_r,u_j,u_j^{-1}]\to \C[\abgr_1]\sim \C[\abgr,\epsilon]\,, \ \
t_r\mapsto \xi_r\,,\ \ u_j\mapsto \xi_j
$$
Using \eqref{map}, we know  that for each $n$, $\theta_n(k)\in
\Z[\abgr_1]$. Now for $k\in \C[t_r,u_j,u_j^{-1}]$ one can compute the
coefficients $b_I$ of the polynomial as the Fourier coefficients
$$
b_I=(2\pi)^{-d}\int_{(S^1)^d}\,k(e^{i\alpha_1},\ldots,e^{i\alpha_d})\,e^{-i
I.\alpha}\prod d\alpha_j
$$
and hence as the limit
$$
b_I=\lim_{n\to \infty} n^{-d}\,\sum_{1\leq a_j\leq n} k(e^{2\pi
i\frac{a_1}{n}},\ldots, e^{2\pi i\frac{a_d}{n}})\,e^{-i I.\alpha}\,, \ \
\alpha_j=2\pi \frac{a_j}{n}
$$
When $k(x)=\prod t_r^{m_r}\prod u_j^{m_j}$ is a monomial, the sum $$
\sum_{1\leq a_j\leq n} k(e^{2\pi i\frac{a_1}{n}},\ldots, e^{2\pi
i\frac{a_d}{n}})\,e^{-i I.\alpha}$$ is either zero or $n^d$ and the latter case
only happens if all the components of the multi index $m-I$ are divisible by
$n$. Thus $$ n^{-d}\,\sum_{1\leq a_j\leq n} k(e^{2\pi i\frac{a_1}{n}},\ldots,
e^{2\pi i\frac{a_d}{n}})\,e^{-i I.\alpha} $$ only depends on $\theta_n(k)$ and
is a relative integer if $\theta_n(k)\in \Z[\abgr_1]$. It follows that all the
$b_I$ are in $\Z$ and hence $k\in \cO(\Omega)$.\endproof

\subsection{The distinction between $G_k$ and  $\mathfrak G(k)$}\label{dist}

Let $\mathfrak G$ be the Chevalley group scheme associated to a root system as
in \S \ref{Cschemes} and let $k$ be a field.

The subgroups $\mathfrak X_r$  generate in the group $\mathfrak G(k)$ of points
that are rational over $k$ a subgroup $G_k$ which is the {\em
commutator subgroup} of the group $\mathfrak G(k)$.
 The subgroup $G_k\subset \mathfrak G(k)$   is often called a Chevalley group over $k$
 and is not in general an algebraic group.
 If $\mathfrak G$ is the universal
Chevalley group, then one knows that $G_k=\mathfrak G(k)$, so that the distinction
between the commutator subgroup $G_k$ and the group $\mathfrak G(k)$ is irrelevant.

In the construction pursued in this paper of the gadget associated to a
Chevalley group, one can take into account this subtlety between $G_k$ and
$\mathfrak G(k)$ by constructing the following sub-gadget. Let $(\tilde
L,\tilde\Phi, n_{\tilde r})$ be the simply connected root system associated to
$(L,\Phi, n_r)$ (\cf~\S~\ref{rootsys}) and let $\varphi:L \to \tilde L$ be the
morphism connecting the two roots systems as follows
\begin{equation}\label{cov}
    \varphi(\Phi)=\tilde\Phi\,, \ \ n_r=n_{\varphi(r)}\circ \varphi\qqq r\in
    \Phi\,.
\end{equation}
 One
simply replaces the term $\Hom(L,\abgr)$ in the construction of the functor
$\underline G(\abgr)$ (\cf~Definition~\ref{defG}) by the following subgroup
\begin{equation}\label{subgad}
\{\chi \in \Hom(L,\abgr)\,|\; \exists \chi'\in \Hom(\tilde L,\abgr)\,, \
\chi=\chi'\circ \varphi\}
\end{equation}
which is the range of the restriction map from $\Hom(\tilde L,\abgr)$ to
$\Hom(L,\abgr)$.

Unlike for the group $\mathfrak G(k)$ the function $N(q)$ that counts the number of
points of $G_k$ for $k=\F_q$ is not in general a polynomial function of $q$.

\section{Schemes over $\F_1$}\label{finrem}

Except for the extra structure given by the grading, the definition of affine variety over $\F_1$ that has been used in this paper is identical to the notion proposed by C. Soul\'e in \cite{Soule} and the replacement of the category $\cR$ by the category $\cF_{ab}$ of finite abelian groups was his suggestion. Thus, Theorem \ref{mainthm} is a solution of the problem formulated in \cite{Soule} (\S~5.4) about Chevalley group schemes.

 Even though this early notion played an important role to get the theory started, it is also too loose inasmuch as the only requirement one imposes is to have enough points with cyclotomic coordinates, but there is no control on the exact size of this set.

In the following, we shall explain that our construction  shows how to strengthen considerably the conditions required  in \cite{Soule} on a  variety over $\F_1$, leading also to the definition of a natural notion of scheme over $\F_1$ which reconciles the original point of view of Soul\'e with the one developed by Deitmar in \cite{deit} following Kurokawa, Ochiai, and Wakayama \cite{KOW} and the log-geometry of monoids of K.~Kato \cite{Kato}.

First, we point out that the category $\cF_{ab}$ of finite abelian groups that has been used in this paper in the definition of the  functor $\underline X$ can be replaced by
the larger category $\cM_{ab}$ of abelian monoids. An abelian monoid is a commutative semi-group with a unit $1$ and a $0$-element. Morphisms in this category send $1$ to $1$ (and $0$ to $0$).

Moreover, the construction of the natural transformation $e_{G}$ can be extended to
yield, for any commutative ring $A$ and for the monoid $M=A$, a map from the
set $\underline
  G(M)$ to the group ${\mathfrak
  G}(A)$. If $A$ is a field, the resulting map is a bijection. More precisely we have the following result

\begin{thm} \label{mainthmbis}  Let ${\mathfrak
  G}$ be the  scheme over $\Z$ associated to a Chevalley group $G$.
\begin{itemize}
  \item The construction \eqref{unG} of the functor $\underline G$ extends to
the category $\cM_{ab}^{(2)}$
 of pairs $(M,\epsilon)$ made by an abelian monoid $M$ and an element $\epsilon$ of square one.
  \item The construction \eqref{egnat} of the natural transformation $e_G$ extends to arbitrary commutative rings $A$  to yield a map
  \begin{equation}\label{mapA}
   e_{
  G,A} :\;\underline G(A,-1) \to {\mathfrak
  G}(A)
 \end{equation}
  \item When $A$ is a field the map $ e_{
  G,A}$ is a bijection.
\end{itemize}
\end{thm}

\proof Let $(M,\epsilon)$ be an object of $\cM_{ab}^{(2)}$.  In  \cite{Tits1}, \S~4.3,  the construction of $\cN_{\abgr,\epsilon}(L,\Phi)$ is carried out for
abelian groups $\abgr$ and it applies in particular to the multiplicative group $M^*$ of a
given monoid $M$. The functor $\underline{G}: \cM_{ab}^{(2)}\to \mathcal Ens$ is the product
    \begin{equation}\label{funcdefn}
    \underline{G}(M,\epsilon)=M^{\Phi^+} \times \coprod_{w\in
    W}(p^{-1}(w)\times M^{\Phi_w})
    \end{equation}
where  $p$ is the projection
$\cN_{M^*,\epsilon}(L,\Phi)\stackrel{p}{\longrightarrow} W$.

One defines the ring $\beta(M,\epsilon)=\Z[M,\epsilon]$ in such a way that the following adjunction relation holds for any commutative ring $A$,
\begin{equation}\label{adjadj}
    \Hom(\beta(M,\epsilon),A)\cong\Hom((M,\epsilon),\beta^*(A))
\end{equation}
where $\beta^*(A)=(A,-1)$ is the object of $\cM_{ab}^{(2)}$ given by the ring $A$ viewed as a
(multiplicative) monoid and the element $-1\in A$. One has
\begin{equation}\label{constadj}
\beta(M,\epsilon)=\Z[M,\epsilon]=\Z[M]/J\,, \ \ J=(1+\epsilon)\Z[M]
\end{equation}
Note that this ring is not always flat over $\Z$ (\eg when $\epsilon=1$) but  the natural morphism $M\to \Z[M,\epsilon]$ is always an injection.

Note also that in \eqref{mapA} the natural transformation sends $\underline G\circ \beta^*$ to
${\mathfrak
  G}$ rather than the way it was formulated above in this paper \ie
$$
 e^{G}:\underline G\to {\mathfrak
  G}\circ\beta,\qquad e^{G}(M,\epsilon):\underline{G}(M,\epsilon) \to {\mathfrak
  G}(\Z[M,\epsilon]).
  $$
  However, the adjunction \eqref{adjadj} shows that these two descriptions are equivalent. Now, we define the natural transformation $ e_{
  G,A}$.
One considers (for each $w\in W$) the maps \eqref{canon1}
 $$
   t=(t_r)_{r\in \Phi^+}\mapsto \psi(t)= \prod_{r\in \Phi^+} x_r (t_r)\in \cU(A)\subset {\mathfrak
  G}(A)
  $$
  and \eqref{canon3}
  $$
    t=(t_r)_{r\in \Phi_w}\mapsto \psi_w(t)=  \prod_{r\in \Phi_w} x_r (t_r)\in \cU_w(A)\subset {\mathfrak
  G}(A).
$$
  For the normalizer one uses Tits' isomorphism (Theorem \ref{titshelp}):
$$
    e_{\cN}\,:\,\cN_{A^*,-1}(L,\Phi)\to \cN(A).
 $$
 For all $n\in p^{-1}(w)$, $a\in M^{\Phi^+}$ and $b\in M^{\Phi_w}$, one lets
 $$
   e_{
  G,A}(a,n,b)=\psi(a)\,e_{\cN}(n)\, \psi_w(b).
  $$
  This definition describes the natural transformation
  $$
  e_{
  G,A} :\;\underline G(A,-1) \to {\mathfrak
  G}(A).
  $$

The last statement follows from Bruhat's Theorem in the form of Theorem
\ref{chevdec}.\endproof

  \smallskip

  Thus Theorem \ref{mainthmbis} suggests to strengthen the conditions imposed
  on a variety over $\F_{1}$ and to define a scheme over $\F_1$ as follows
  (using the adjoint pair of functors $\beta(M)=\Z[M]$ and $\beta^*(A)=A$)

  \begin{defn}\label{defnfonesch} An $\F_1$-scheme is given by a covariant functor $\underline X$ from $\cM_{ab}$ to the category of sets, a $\Z$-scheme $X$ and a natural transformation from
  $\underline X\circ \beta^*$ to $X$ such that
\begin{itemize}
    \item $\underline X$ is locally representable.
  \item The natural transformation  is bijective if $A=\K$ is a field.
\end{itemize}
\end{defn}
  For arbitrary commutative rings the natural transformation
  yields a map
  \begin{equation}\label{strongcond}
    e_{X,A}\,:\; \underline X(A)\to X(A)
  \end{equation}
which is bijective when $A$ is a field.

A nice feature of this definition is
that it ensures that the counting of points gives the correct answer. Indeed,
the above conditions ensure that the number of points over $\F_{1^n}$ which is
given by the cardinality of $\underline X(\abgr)$ for\footnote{One   adjoins a zero element
to get a monoid in the above sense} $\abgr=\Z/n\Z$, agrees
with the cardinality of $X(\F_q)$ when $n=q-1$ and $q$ is a prime power.

In the forthcoming paper \cite{ak1}, we shall develop the general theory of $\F_1$-schemes  and show in particular that under natural finiteness conditions the restriction of the functor $\underline X$
to the full subcategory $\cF_{ab}$ of $\cM_{ab}$ of finite abelian groups is automatically a functor to finite graded sets. Moreover under a torsion free hypothesis the number of points over $\F_{1^n}$ is a polynomial $P(n)$, where $P(x)$ has positive integral coefficients.
 Unlike the theory developed in \cite{deit}, our conditions do not  imply that the varieties under study are necessarily toric.

\smallskip
 What we have shown in this paper is that Chevalley group schemes yield
 varieties over $\F_{1^2}$, but {\em not} that the group operation $\mu$ can be
 defined over $\F_{1^2}$. In fact, it is only the terms of lowest degree (equal
 to the rank of $G$) that yield a group, namely the group $\cN_{\abgr,\epsilon}(L,\Phi)$
 of J. Tits. The structure of the terms of higher order is more mysterious.

\end{document}